\documentclass[12pt]{article}
\usepackage[top=1in, bottom=1in, left=1in, right=1in]{geometry}

\usepackage{xcolor}
\usepackage{geometry}
\usepackage{amssymb,bm}
\usepackage{amsmath}
\usepackage{amsthm}
\usepackage{graphicx}
\usepackage{hyperref}
\usepackage{caption}
\usepackage{float} 
\hypersetup{pdfborder=0 0 0}
\usepackage{calrsfs}
\usepackage{booktabs}
\usepackage{physics}

\newtheorem{theorem}{{\sc Theorem}}[section]

\newtheorem{lemma}[theorem]{{\sc Lemma}}

\newtheorem{remark}[theorem]{Remark}
\newtheorem{definition}[theorem]{Definition}

\newcommand\restr[2]{{ \left.\kern-\nulldelimiterspace #1 \vphantom{\big|} \right|_{#2}}}
\newcommand{\RR}{\mathbb{R}}

\newcommand{\CC}{\mathbb{C}}

\newcommand{\n}{\noindent}
\newcommand{\comment}[1]{}

\DeclareMathAlphabet{\pazocal}{OMS}{zplm}{m}{n}

\newcommand{\BS}{\mathbb{S}}

\usepackage{bbm}
\usepackage{subcaption}

\title{Scattering of plane waves and a modified Pompeiu criterion}
\author{
Narek Hovsepyan\footnote{School of Mathematical and Statistical Sciences, University of Texas Rio Grande Valley, Edinburg, TX, USA (narek.hovsepyan@utrgv.edu)} \quad and \quad Michael S. Vogelius\footnote{Department of Mathematics, Rutgers University, New Brunswick, NJ, USA (vogelius@math.rutgers.edu)}}
\date{}

\begin{document}
\maketitle

\begin{abstract}

In this paper we prove that plane waves always scatter from bounded, strictly convex planar inhomogeneities with constant material parameters. Our approach involves the introduction of a modified Pompeiu property, whose failure corresponds to the vanishing of the Fourier transform of the indicator function of the inhomogeneity on a complexified circle not centered at the origin. We then show that such vanishing is impossible for bounded, strictly convex planar sets.  
\end{abstract}

\section{Introduction} \label{Intro}

This paper has two separate but very related goals: (1) to show that any strictly convex inhomogeneity with constant index of refraction scatters any plane wave non-trivially, and in doing so (2) to elucidate the similarities, as well as the differences, when compared to the well known Pompeiu/Schiffer problems. The scattering problem we consider is described by the two dimensional Helmholtz equation

\begin{equation*}
\nabla \cdot (\gamma(\bm{x}) \nabla U) +k^2 n(\bm{x}) U=0 \qquad \text{in} \ \mathbb{R}^2,
\end{equation*}

\n where

\begin{equation*}
\gamma(\bm{x}) =
\begin{cases}
a, \qquad \bm{x} \in D
\\
1, \qquad \text{otherwise},
\end{cases}
\qquad \text{and} \qquad
n(\bm{x}) =
\begin{cases}
q, \qquad \bm{x} \in D
\\
1, \qquad \text{otherwise}.
\end{cases}
\end{equation*}

\n Here $D$ is the bounded domain representing the inhomogenity, and $a$ and $q$ are positive constants different from $(1,1)$. If $u^{\text{in}}$ denotes the incident wave (a solution to $\Delta u^{\text{in}}+k^2u^{\text{in}}=0$ in all of $\mathbb{R}^2$) then the total field, $U$, decomposes as $U=u^{\text{in}}+u^{\text{sc}}$, where the scattered part satisfies the so called Sommerfeld Radiation Condition. We consider incident plane waves, that is waves of the form $u^{\text{in}}(\bm{x})=e^{ik\bm{x}\cdot \bm{\eta}}$ for some unit vector $\bm{\eta} \in \RR^2$. As one of the main results of this paper we show that if $D$ is strictly convex (or just convex if $a=1$) then $u^{\text{sc}}$ does not vanish identically outside $D$. This result significantly extends the results obtained in \cite{HV2} and \cite{VX}. The path towards this result is to examine the Fourier Transform of the indicator function of the set $D$ and in particular the structure of the zeroes of this transform. We first show that the vanishing of $u^{\text{sc}}$ outside $D$ (nonscattering) implies the vanishing of this Fourier Transform on the complexified circle 
$$
S_\kappa(k\bm{\eta})=\{ \bm{z}\in \mathbb{C}^2~:~ (\bm{z}- k\bm{\eta}) \cdot (\bm{z}- k\bm{\eta})=\kappa^2 \} ~,
$$
where $\cdot$ denotes the inner product of $\RR^2$ and $\kappa$ denotes the wavenumber in the medium $D$ given by 

\begin{equation} \label{kappa}
\kappa = k \sqrt{\frac{q}{a}}.    
\end{equation}

\n This is completely analogous to the conclusion which one arrives at for sets that do not have the Pompeiu Property, except that the complexified circle would be centered at the origin. The other main result of this paper then shows that for a strictly convex domain such vanishing of the Fourier Transform on $S_\kappa(k\bm{\eta})$ is impossible. The reason this is easier to do for our scattering problem than for the Pompeiu/Schiffer problems is that the vanishing of the Fourier Transform on $S_\kappa(k\bm{\eta})$ implies the vanishing on $S_\kappa(-k\bm{\eta})$ as well. This fact (at least for $\kappa > k$) leads to the existence of a zero of multiplicity 2 for the Fourier Transform. Once this is achieved, techniques along the lines of those used in \cite{BK} and \cite{HV2} lead to a contradiction. For the case $\kappa \leq k$ it was already proven in \cite{HV2} that any bounded, simply connected, planar Lipschitz domain scatters any plane wave.

\section{Preliminaries} \label{SECT Prelim}
\setcounter{equation}{0}

Throughout, we assume that $D \subset \RR^2$ is a nonempty, bounded and simply connected domain. For such a domain, consider the Fourier transform of its characteristic function

\begin{equation*}
\hat{\chi}_D(\bm{z}) = \int_D e^{i \bm{x} \cdot \bm{z}} d\bm{x}.
\end{equation*}

\n We introduce the following notation.

\begin{definition}
\n Given $\bm{c} \in \RR^2$ and $\kappa>0$, the complexified circle of radius $\kappa$ centered at $\bm{c}$ is defined by\footnote{We note that the centers of the complexified circles introduced here are always real, as in the scattering application $\bm{c} = k \bm{\eta} \in \RR^2$.}

\begin{equation*}
S_\kappa(\bm{c}) = \left\{ \bm{z} \in \CC^2 ~:~ (\bm{z} - \bm{c}) \cdot (\bm{z}-\bm{c}) = \kappa^2 \right\},
\end{equation*}

\n where $\cdot$ denotes the inner product of $\RR^2$.
\end{definition}

Motivated by the study of scattering of plane waves from the inhomogeneity $D$ as discussed in the introduction, we pose the following question: can there exist $\bm{c} \in \RR^2$ and $\kappa>0$, such that

\begin{equation} \label{pomp c}
\hat{\chi}_D = 0 \qquad \text{on} \ S_\kappa(\bm{c}).
\end{equation}

\n This can be viewed as a modified version of the failure of the Pompeiu property, with the complexified circle allowed to be off-centered when $\bm{c} \neq \bm{0}$. Indeed, $D$ fails to have the Pompeiu property \cite{BST73} if and only if there exists $\kappa > 0$, such that

\begin{equation} \label{pomp}
\hat{\chi}_D = 0 \qquad \text{on} \ S_\kappa(\bm{0}).
\end{equation}

\n It is well known that the disk fails to have the Pompeiu property. The celebrated Pompeiu (or Schiffer) conjecture asserted that the disk is the only exceptional simply connected domain that fails to have the Pompeiu property. Equivalently, for any noncircular simply connected domain $D$, $\hat{\chi}_D$ cannot vanish identically on a complexified circle centered at the origin. The recent works \cite{CS, conject} disproved this conjecture by constructing noncircular simply connected planar domains $D$ that fail to have the Pompeiu property. These counterexamples, however, are nonconvex and it remains an open question whether the Pompeiu conjecture holds within the class of strictly convex domains.

As part of this work, we study the modified question \eqref{pomp c} for $\bm{c} \neq \bm{0}$ and in that case we show that every bounded strictly convex planar $D$, including the disk, has the off-centered Pompeiu property. More precisely, we prove that when $\bm{c} \neq \bm{0}$, the identity \eqref{pomp c} cannot hold. We will see that the ``broken symmetry" arising from the complexified circle being off-centered has a strong consequence: it imposes an additional constraint that allows us to establish the result. This constraint degenerates in the centered case $\bm{c}=\bm{0}$, corresponding to the classical Pompeiu property, and therefore our approach does not apply in that setting. As a consequence of this result, we show that plane waves always scatter from any bounded strictly convex inhomogeneity $D$ with constant material parameters $(a,q) \neq (1,1)$. In \cite{HV2}, we established this always-scattering result for elongated strictly convex $D$ in the case $a>k/\kappa$, building on the ideas of \cite{BK}, where it was shown that any elongated strictly convex $D$ has the Pompeiu property. In the present work, we remove the elongatedness assumption as well as the restriction on $a$.

We return to the scattering problem. Consider an incident plane wave with propagation direction given by the unit vector $\bm{\eta}$ and wavenumber $k>0$:

\begin{equation} \label{u^in}
u^{\text{in}}(\bm{x}) = e^{ik \bm{\eta} \cdot \bm{x}}, \qquad \qquad \bm{\eta} \in \BS = \{\bm{\eta} \in \RR^2: |\bm{\eta}| = 1\}.
\end{equation}  

\n Let $D \subset \RR^2$ be a bounded, simply connected Lipschitz domain with constant material parameters $a,q>0$ with $(a,q) \neq (1,1)$. We say that $u^{\text{in}}$ is nonscattering if and only if the scattered wave generated in the presence of the inhomogeneity $D$ vanishes identically outside $D$. Equivalently, this means that the following problem admits a solution in $D$:

\begin{equation} \label{u system}
\begin{cases}
\Delta u^{\text{sc}} + \kappa^2 u^{\text{sc}} = (k^2-\kappa^2)u^{\text{in}} \quad
& \text{in} \ D,
\\
u^{\text{sc}} = 0 & \text{on} \ \partial D,
\\
a \, \partial_{\bm{\nu}} \left( u^{\text{sc}} + u^{\text{in}}\right) = \partial_{\bm{\nu}} u^{\text{in}}  & \text{on} \ \partial D,
\end{cases}
\end{equation}

\n where $\kappa$ is given by \eqref{kappa} and $\bm{\nu}$ denotes the unit outer normal to $\partial D$. Multiplying the PDE in \eqref{u system} by the test function $\phi(\bm{x}) = e^{i\bm{x}\cdot\bm{\zeta}}$ and using Green's identity, we obtain \cite{HV2}

\begin{equation*}
k^2 (1-q) \int_D u^{\text{in}}(\bm{x}) e^{i\bm{x}\cdot\bm{\zeta}} d\bm{x} - i(1-a) \int_D \nabla u^{\text{in}}(\bm{x}) \cdot \bm{\zeta} e^{i\bm{x}\cdot\bm{\zeta}} d\bm{x} = 0 \qquad \quad \forall \, \bm{\zeta} \in S_\kappa(\bm{0}).
\end{equation*}

\n Using the plane wave form \eqref{u^in} of the incident field and setting $\bm{z} = k\bm{\eta} + \bm{\zeta}$, we obtain:

\begin{equation}
C(\bm{z}) \hat{\chi}_D(\bm{z}) = 0 \qquad \qquad \forall \ \bm{z} \in S_{\kappa}(k \bm{\eta}),
\end{equation}

\n where

\begin{equation}
C(\bm{z}) = a-q + \frac{1-a}{k} \bm{z} \cdot \bm{\eta}.
\end{equation}

\n When $a=1$, $C(\bm{z})$ is never zero and consequently $\hat{\chi}_D$ vanishes identically on $S_{\kappa}(k \bm{\eta})$. If $a \neq 1$, then $C(\bm{z}) = 0$ if and only if

\begin{equation} \label{z eta = alpha}
\bm{z} \cdot \bm{\eta} = k \frac{q-a}{1-a}.
\end{equation}

\n There are at most two points on the complexified circle $S_{\kappa}(k \bm{\eta})$ that satisfy \eqref{z eta = alpha}. Therefore, $\hat{\chi}_D$ vanishes on $S_{\kappa}(k \bm{\eta})$ except possibly at these two points. However, since $\hat{\chi}_D$ is continuous (in fact, entire in $\CC^2$), it must vanish on the entire complexified circle. Thus, we obtain:

\begin{lemma} \label{LEM nonsc to pomp}
Assume the setup introduced above and suppose that the plane wave \eqref{u^in} is nonscattering for the inhomogeneity $D$. Then

\begin{equation} \label{chi_D vanish}
\hat{\chi}_D  = 0 \qquad \text{on} \ S_\kappa(k\bm{\eta}).
\end{equation}

\end{lemma}

\begin{remark} \normalfont

In the case $a=1$, the condition \eqref{chi_D vanish} is, in fact, equivalent to the existence of a solution to \eqref{u system} with $u^{\text{in}}(\bm{x}) = e^{ik \bm{\eta} \cdot \bm{x}}$. This follows from the fact that when $a=1$, the function $u^{\text{sc}}$ in \eqref{u system} has both zero Dirichlet- and zero Neumann trace. One can then proceed using the same arguments as in \cite{w1}, see also the appendix in \cite{HV1}. However, we do not use this reverse implication here.  

\end{remark}

\section{Main results} 
\setcounter{equation}{0}

\begin{theorem} \label{THM main}
Let $D \subset \RR^2$ be a nonempty bounded domain, and let $\bm{c} \in \RR^2$ and $\kappa>0$. Then the following statements hold:

\begin{enumerate}
\item[(i)] If $|\bm{c}| \geq \kappa$, then $\hat{\chi}_D$ cannot vanish identically on $S_\kappa (\bm{c})$.

\item[(ii)] If $0<|\bm{c}| < \kappa$, then $\hat{\chi}_D$ cannot vanish identically on $S_\kappa (\bm{c})$, provided $D$ is strictly convex. 
\end{enumerate}

\end{theorem}

\n Thus, for strictly convex $D$, combining both parts of the above theorem shows that $D$ has the off-centered Pompeiu property, i.e. $\hat{\chi}_D$ cannot vanish identically on any complexified circle $S_\kappa (\bm{c})$ with $\bm{c} \in \RR^2 \backslash \{\bm{0}\}$ and $\kappa>0$. In view of Lemma~\ref{LEM nonsc to pomp}, we therefore obtain

\begin{theorem}
Let $D \subset \RR^2$ be a nonempty, bounded, strictly convex inhomogeneity with constant material parameters $a,q>0$ satisfying $(a,q) \neq (1,1)$. Then every incident plane wave scatters from $D$. More precisely, \eqref{u^in} scatters for every wavenumber $k>0$ and every unit vector $\bm{\eta} \in \RR^2$. 
\end{theorem}

In the case $a=1$, the above always-scattering result holds for convex inhomogeneities, rather than only for strictly convex ones. This follows due to the regularity result of \cite{CV}, which in the present setting implies that if there exists a plane wave that does not scatter, then the boundary of $D$ must be real analytic. In particular, if $D$ is convex, real analyticity prevents its boundary from containing any line segments and therefore implies that $D$ is strictly convex. We are then in a position to apply Theorem~\ref{THM main}, which, together with Lemma~\ref{LEM nonsc to pomp} establishes the theorem. In the case $a\neq 1$, the regularity result of \cite{CVX23} holds for $C^{1,\alpha}$ boundaries and applies at boundary points where the propagation direction $\bm{\eta}$ is not tangential to the boundary. The extension to Lipschitz boundaries, established in \cite{KSS24}, applies at boundary points having a neighborhood in which $\bm{\eta}\cdot\bm{\nu}$ has a fixed sign and is uniformly bounded away from zero. In summary, for $a \neq 1$ we cannot automatically conclude from nonscattering of a plane wave that the boundary is real analytic.

Part $(i)$ of Theorem~\ref{THM main} is not new; it was established in our previous work \cite{HV2} (see Theorem 3.1 of that paper). Here, however, it is stated in a slightly different form, namely, in terms of the modified Pompeiu property rather than scattering. Since the proof is quite short, we include it below in a form adapted to the present setting. The key component is part $(ii)$ of Theorem~\ref{THM main}, which follows directly from the two lemmas below.

\begin{lemma} \label{LEM double zero}
Let $D \subset \RR^2$ be a nonempty and bounded domain. Suppose that there exist $\bm{c} \in \RR^2$ and $\kappa > 0$ with $0<|\bm{c}| < \kappa$, such that $\hat{\chi}_D = 0$ on $S_\kappa(\bm{c})$. Then there exists a point $\bm{z}_0 \in \RR^2$, such that

\begin{equation} \label{double zero}
\hat{\chi}_D(\bm{z}_0) = \nabla \hat{\chi}_D(\bm{z}_0) = 0.
\end{equation}

\end{lemma}

\begin{lemma} \label{LEM strict convex}
Let $D \subset \RR^2$ be a nonempty, bounded and strictly convex domain. Then $\hat{\chi}_D$ cannot have a double zero in $\RR^2$, i.e., \eqref{double zero} cannot hold.
\end{lemma}

Regarding Lemma~\ref{LEM double zero}, we remark that the condition $\bm{c} \neq \bm{0}$ is essential for the existence of a double zero. In the case $\bm{c} = \bm{0}$ our construction of $\bm{z}_0$ fails. Regarding Lemma~\ref{LEM strict convex}, we remark that if the strict convexity assumption is dropped, then double zeros may occur. Indeed, consider the square $D = (-1,1)^2$. Then, writing $\bm{z}=(z_1, z_2)$ and similarly for $\bm{x}$, we have

\begin{equation*}
\hat{\chi}_D(\bm{z}) = \int_{-1}^1 e^{iz_1 x_1} dx_1  \int_{-1}^1 e^{iz_2 x_2} dx_2 = f(z_1) f(z_2),
\end{equation*}

\n where $f(x) = 2 \sin x / x$. It is clear that $\hat{\chi}_D$ has a double zero at $\bm{z}_0 = (\pi, \pi)$.

\section{Proofs} 
\setcounter{equation}{0}

\subsection{Part $(i)$ of Theorem~\ref{THM main}}

Assume, to the contrary, that $\hat{\chi}_D$ vanishes identically on $S_\kappa(\bm{c})$ for some $|\bm{c}| \geq \kappa > 0$. We will show that in this case

\begin{equation} \label{intersect}
S_\kappa(\bm{c}) \cap i \RR^2 \neq \emptyset.
\end{equation}

\n This concludes the proof. Indeed, if we denote the intersection point by $i \bm{y}$, where $\bm{y} \in \RR^2$, then

\begin{equation*}
\hat{\chi}_D(i\bm{y}) = \int_D e^{- \bm{x} \cdot \bm{y}} d \bm{x} > 0,
\end{equation*}

\n which gives the desired contradiction. It remains to prove \eqref{intersect}. For $i \bm{y} \in S_\kappa(\bm{c})$ we require  

\begin{equation*}
(i \bm{y} - \bm{c}) \cdot (i \bm{y} - \bm{c}) = \kappa^2,
\end{equation*}

\n or, equivalently,

\begin{equation*}
|\bm{y}|^2 + 2i \bm{y} \cdot \bm{c} = |\bm{c}|^2 - \kappa^2.
\end{equation*}

\n It is clear that such a $\bm{y}$ exists. Indeed, take $\bm{y} \perp \bm{c}$ with $|\bm{y}| = \sqrt{|\bm{c}|^2 - \kappa^2}$.

\subsection{Proof of Lemma~\ref{LEM double zero}}

Recall that, by hypothesis, $\hat{\chi}_D = 0$ on $S_\kappa(\bm{c})$ for some $0<|\bm{c}| < \kappa$. Our goal is to prove that $\hat{\chi}_D$ has a double zero at some point $\bm{z}_0 \in \RR^2$. For this purpose, we will not use the full complexified circle, but only its real counterpart

\begin{equation*}
R_\kappa(\bm{c}) = S_\kappa(\bm{c}) \cap \RR^2.
\end{equation*}

\n The proof is based on symmetry. Since $\chi_D$ is real valued, its Fourier transform satisfies $\overline{\hat{\chi}_D(\bm{z})} = \hat{\chi}_D(- \overline{\bm{z}})$. Consequently, the hypothesis implies that 

\begin{equation*}
\hat{\chi}_D(- \overline{\bm{z}})=0 \qquad \qquad \forall \, \bm{z} \in R_\kappa(\bm{c})
\end{equation*}

\n or, equivalently,

\begin{equation*}
\hat{\chi}_D(\bm{z})=0, \qquad \qquad \forall \, \bm{z} \in R_\kappa(-\bm{c}).
\end{equation*}

\n Thus, $\hat{\chi}_D$ vanishes on two circles:

\begin{equation*}
\hat{\chi}_D(\bm{z})=0, \qquad \qquad \forall \, \bm{z} \in R_\kappa(\bm{c}) \cup R_\kappa(-\bm{c})
\end{equation*}

\n Since $\kappa> |\bm{c}|$, these two circles intersect and we may take $\bm{z}_0$ to be either of the two intersection points. Clearly, $\hat{\chi}_D(\bm{z}_0) = 0$. To show that the gradient also vanishes at $\bm{z}_0$, we use the fact that the two circles intersect nontangentially. We remark that in the classical Pompeiu case, corresponding to $\bm{c} = \bm{0}$, the two circles coincide and hence there is no nontangential intersection.

More precisely, let $\bm{t}_+$ and $\bm{t}_-$ denote tangent vectors to the circles $R_\kappa(\bm{c})$ and $R_\kappa(-\bm{c})$, respectively, at their common point $\bm{z}_0$. Since $\hat{\chi}_D$ vanishes identically on $R_\kappa(\pm\bm{c})$, its tangential derivatives along the two circles also vanish at $\bm{z}_0$:

\begin{equation*}
\nabla \hat{\chi}_D(\bm{z}_0) \cdot \bm{t}_\pm = 0.
\end{equation*}

\n The vectors $\bm{t}_+$ and $\bm{t}_-$ are linearly independent, and therefore $\nabla \hat{\chi}_D(\bm{z}_0) = 0$.

\subsection{Proof of Lemma~\ref{LEM strict convex}}

For the sake of contradiction, suppose that $\hat{\chi}_D$ has a double zero at some $\bm{z}_0 \in \RR^2$. Clearly, $\bm{z}_0 \neq 0$, since otherwise $D = \emptyset$. We may therefore rotate the coordinate axes so that $\bm{z}_0 = (x_0, 0)$ with $x_0>0$. Now consider the translated domain $D' = \bm{v} + D$, where $\bm{v} \in \RR^2$. By the standard properties of the Fourier transform, using a change of variables and the product rule, we obtain

\begin{equation*}
\hat{\chi}_{D'}(\bm{z}) = e^{i \bm{v} \cdot \bm{z}} \hat{\chi}_{D}(\bm{z}), \qquad \qquad
\nabla \hat{\chi}_{D'}(\bm{z}) = e^{i \bm{v} \cdot \bm{z}} \left[ \nabla \hat{\chi}_{D}(\bm{z}) + i \bm{v} \hat{\chi}_{D}(\bm{z}) \right].
\end{equation*}

\n In particular, the double zero $\bm{z}_0$ of $\hat{\chi}_D$ is preserved under translations of $D$. We now consider the horizontal projection of $D$ onto $x$-axis. In view of this translation invariance, we may translate the coordinate system so that this projection is the interval $[0, w]$, where $w>0$ denotes the horizontal width of $D$. Thus, we may write

\begin{equation*}
\overline{D} = \left\{ (t, s): 0 \leq t \leq w, \quad L_1(t) \leq s \leq L_2(t) \right\},
\end{equation*}

\n where $L_1$ and $L_2$ denote the lower and upper portions of the boundary, respectively. Since $D$ is strictly convex, $L_1$ is strictly convex, while $L_2$ is strictly concave. Consequently, the difference $L = L_2 - L_1$ is strictly concave. In particular, $L'$ exists almost everywhere and is integrable on $(0, w)$.
 
We will show that the existence of a double zero of $\hat{\chi}_D$, together with the strict convexity of $D$, implies that there exists $r>0$ such that

\begin{equation} \label{f sine = 0}
\int_0^r f(t) \sin t dt = 0,
\end{equation}

\n where $f$ is a positive, strictly decreasing and integrable function on $(0,r)$, expressed in terms of $L$ (see below). We will then obtain a contradiction by showing that the above integral cannot vanish. We state this fact as a lemma, whose proof is postponed to the next subsection.

\begin{lemma} \label{LEM Polya}
Let $r>0$ and $f$ be a positive, strictly decreasing and integrable function on $(0, r)$, then

\begin{equation}
\int_0^r f(t) \sin t dt > 0
\end{equation}

\end{lemma}

\begin{remark}
\normalfont
In \cite{HV2} (see also \cite{BK}), we already encountered this type of sine integral and obtained a contradiction to \eqref{f sine = 0} under appropriate conditions on $r$ (see, for example, the proof of Theorem 3.7 in \cite{HV2}). There, however, the function playing the role of $f$ was $L'$, which is strictly decreasing but does not possess the additional positivity enjoyed by $f$ here. This additional positivity will allow us to obtain the desired contradiction without imposing any conditions on $r$. We also note that the above lemma is proved in \cite{polya}. The study of such sine integrals and their zeros goes back to P\'olya.
\end{remark}

Let us show how to obtain \eqref{f sine = 0} based on the existence of the double zero $\bm{z}_0=(x_0, 0)$. Note that, with $\bm{z} = (z_1, z_2)$, we have

\begin{equation*}
\hat{\chi}_D(\bm{z}) = \int_0^w e^{iz_1 t} \int_{L_1(t)}^{L_2(t)} e^{i z_2 s} ds dt, 
\qquad \qquad
\partial_{z_1} \hat{\chi}_D(\bm{z}) = i \int_0^w t e^{iz_1 t} \int_{L_1(t)}^{L_2(t)} e^{i z_2 s} ds dt.
\end{equation*}

\n We now evaluate these functions on the horizontal axis $z_2=0$. To this end, set

\begin{equation*}
F(x) = \hat{\chi}_D((x,0)),
\end{equation*}

\n then

\begin{equation*}
F(x) = \int_0^w L(t) e^{ixt} dt, \qquad \qquad F'(x) = i \int_0^w t L(t) e^{ixt} dt.
\end{equation*}

\n By the double zero hypothesis

\begin{equation} \label{F and F' zero}
F(x_0) = F'(x_0) = 0.
\end{equation}

\n Integrating by parts in the integral defining $F(x)$ and using $L(0)=L(w)=0$, we obtain

\begin{equation*}
F(x) = - \int_0^w \left[ L(t) e^{ixt} \right]' t dt = - \int_0^w t L'(t)e^{ixt} dt - x F'(x).
\end{equation*}

\n Therefore, \eqref{F and F' zero} implies that

\begin{equation} \label{1}
\int_0^w L(t) e^{ix_0t} dt = 0 \qquad \text{and} \qquad \int_0^w t L'(t) e^{ix_0t} dt = 0.
\end{equation}

\n Consider the integrable function

\begin{equation} \label{g}
g(t) = L(t) - tL'(t),
\end{equation}

\n defined almost everywhere in $(0,w)$. Below, we will show that it admits a representative, still denoted by $g$, which is positive and strictly increasing on $(0,w)$. Subtracting the two equations in \eqref{1} then gives

\begin{equation*}
\int_0^w g(t) e^{ix_0t} dt = 0.
\end{equation*}

\n To obtain the desired equation \eqref{f sine = 0}, we replace $g$ by a positive, strictly decreasing function. This is achieved by setting $\tilde{f}(t) = g(w-t)$. A change of variables in the above integral and cancellation of the phase factor $e^{i x_0 w}$ now gives

\begin{equation*}
\int_0^w \tilde{f}(t) e^{-ix_0t} dt = 0.
\end{equation*}

\n Taking imaginary parts we arrive at

\begin{equation*}
\int_0^w \tilde{f}(t) \sin (x_0 t) dt = 0,
\end{equation*}

\n where $\tilde{f}$ is positive and strictly decreasing. Changing variables once more, we obtain

\begin{equation*}
\int_0^{r} f(t) \sin t dt = 0, \qquad \qquad r = w x_0,
\end{equation*}

\n with $f(t) = \tilde{f}(t/x_0)$. The latter is precisely \eqref{f sine = 0}.

Thus, it only remains to establish that the function $g$ defined by \eqref{g} has a positive and strictly increasing representative. Since $L$ is concave, it has one-sided derivatives at every point of $(0,w)$. Let $L'_+$ denote its right derivative, and consider the following representative of the function \eqref{g}:

\begin{equation*}
g(t) = L(t) - t L'_+(t), \qquad \qquad t \in (0, w).
\end{equation*}

\n Let us show that this function has the desired properties.

\vspace{.05in}

$\bullet$ $g$ is positive. To prove this, we use the fact that since $L$ is strictly concave on $[0,w]$, its right derivative $L'_+$ is strictly decreasing on $(0,w)$ \cite{rock}. Fix any $t \in (0,w)$. Since $L'$ agrees with $L'_+$ almost everywhere, we may write

\begin{equation} \label{L(t)-L(0)}
L(t) - L(0) = \int_0^t L_+'(x) dx > \int_0^t L_+'(t) dx = t L_+'(t).
\end{equation}

\n Using $L(0) = 0$, we conclude that $g(t)>0$.

\vspace{.05in}

$\bullet$ $g$ is strictly increasing. To prove this, take $0<a<b<w$. First, note that, by the same argument used to obtain \eqref{L(t)-L(0)}, we have

\begin{equation*}
L(b) - L(a) > (b-a) L'_+(b),
\end{equation*}

\n but then

\begin{equation*}
g(b) - g(a) = L(b)-L(a) + a L'_+(a) - b L'_+(b) > a \left[ L'_+(a) - L'_+(b) \right] > 0 .
\end{equation*}

\subsubsection{Proof of Lemma~\ref{LEM Polya}}

For the moment, suppose that $r=2\pi$. Then, as in \cite{HV2}, we split

\begin{equation} \label{int pos}
\int_0^{2 \pi} f(t) \sin t dt = \int_0^\pi + \int_\pi^{2\pi} = \int_0^\pi \left[ f(t) - f(\pi + t) \right] \sin t dt  > 0.
\end{equation}

\n Note that this inequality uses only the strict monotonicity of $f$ and its positivity was not used. For general $r$, we split $(0,r)$ into intervals of length $\pi$, pair consecutive integrals as above and use the positivity of $f$ to handle the remaining ``endpoint" interval. We assume that $r>\pi$, since otherwise the claim is immediate from the positivity of the integrand. Dividing $r$ by $\pi$, we may write

\begin{equation*}
r = m \pi + r_0 \qquad \qquad m \in \mathbb{N}, \quad r_0 \in [0, \pi).
\end{equation*}

\n If $m = 2 n$ is even, then we split 

\begin{equation*}
\int_0^r = \sum_{j=1}^{n} \int_{2\pi (j-1)}^{2\pi j} + \int_{2\pi n}^{2\pi n + r_0}.
\end{equation*}

\n Each integral in the sum is positive by \eqref{int pos}, while the last integral is nonnegative because $f$ is positive and so is $\sin$, since $r_0 \in [0,\pi)$.

If $m=2n+1$ is odd, then 

\begin{equation*}
r = 2\pi n + \pi + r_0
\end{equation*}

\n and we split 

\begin{equation*}
\int_0^r = I_1 + I_2, \qquad \text{where} \qquad I_1 = \int_0^{2\pi n} \quad \text{and} \quad I_2 = \int_{2\pi n}^{2\pi n + \pi} + \int_{2\pi n +\pi}^r.
\end{equation*}

\n Clearly, $I_1 \geq 0$ by the preceding argument for even $m$. It therefore remains to show that $I_2>0$. It suffices to consider the case $n=0$, since shifting the integration variable by $2\pi n$ leaves the sine function unchanged and preserves the positivity and strict monotonicity of $f$. Thus, let $n=0$, so that $r = \pi + r_0 < 2\pi$. Extending $f$ by zero on $(r,2\pi)$, we obtain

\begin{equation*}
I_2 = \int_0^\pi + \int_\pi^r = \int_0^\pi + \int_\pi^{2\pi} = \int_0^\pi \left[ f(t) - f(\pi + t) \right] \sin t dt.
\end{equation*}

\n The last integral is again positive because its integrand is positive. Indeed, if $\pi + t > r$, then $f(t) - f(\pi + t) = f(t) > 0$. On the other hand, if $\pi + t < r$, then both terms are nonzero but $f(t) - f(\pi + t) > 0$ since $f$ is strictly decreasing.

\bibliographystyle{abbrv}
\bibliography{ref}

\end{document}